\documentclass[10pt,reqno]{amsart}
\usepackage{amsmath}
  \usepackage{paralist}
  \usepackage{graphics} %% add this and next lines if pictures should be in esp format
  \usepackage{epsfig} %For pictures: screened artwork should be set up with an 85 or 100 line screen
  \usepackage{graphicx}
  \usepackage{epstopdf}%This is to transfer .eps figure to .pdf figure; please compile your paper using PDFLeTex or PDFTeXify.
  \usepackage[colorlinks=true]{hyperref}
  \usepackage{comment}
  \hypersetup{urlcolor=blue, citecolor=red}

\newtheorem{theorem}{Theorem}[section]

\newtheorem{lemma}[theorem]{Lemma}
\newtheorem{proposition}{Proposition}

\theoremstyle{definition}
\newtheorem{definition}[theorem]{Definition}
\newtheorem{remark}{Remark}

\newcommand{\eps}[1]{{#1}_{\varepsilon}}

\newtheorem{meta-thm}[theorem]{Meta-Theorem}

\newtheorem{algorithm}[theorem]{Algorithm}
\newcommand{\noaverage}[1]{({#1})^0}

\newcommand\beq[1]{ \begin{equation}\label{#1} }
\newcommand{\eeq}{ \end{equation} }

\newcommand\beqa[1]{ \begin{eqnarray} \label{#1}}

\newcommand{\eeqa}{ \end{eqnarray} }

\newcommand{\beqano}{ \begin{eqnarray*} }

\newcommand{\eeqano}{ \end{eqnarray*} }

\newcommand\equ[1]{{\rm (\ref{#1})}}

\def\dist{\operatorname{dist}}

\def\Id{\operatorname{Id}}
\def\Im{\operatorname{Im}}

\def\Range{\operatorname{Range}}
\def\A{{\mathcal A}}
\def\B{{\mathcal B}}

\def\D{{\mathcal D}}
\def\G{{\mathcal G}}
\def\E{{\mathcal E}}
\def\tE{{\tilde  E}}

\def\I{{\mathcal I}}
\def\M{{\mathcal M}}

\def\S{{\mathcal S}}
\def\T{{\mathcal T}}
\def\U{{\mathcal U}}
\def\complex{{\mathbb C}}
\def\integer{{\mathbb Z}}
\def\nat{{\mathbb N}}

\def\real{{\mathbb R}}
\def\torus{{\mathbb T}}

\def\eps{\varepsilon}
\def\th{\theta}

\def\st{s}
\def\un{u}
\title[Whiskered KAM tori of conformally symplectic systems]
{Whiskered KAM tori of conformally symplectic systems}

\author[Renato C. Calleja and Alessandra Celletti and Rafael de la Llave]{}


\email{calleja@mym.iimas.unam.mx}
\email{celletti@mat.uniroma2.it}
\email{rll6@math.gatech.edu}

\thanks{R.C. was partially supported by PAPIIT-DGAPA grant IA102818.
  A.C. was partially supported by GNFM-INdAM and
acknowledges the MIUR Excellence Department Project awarded to the Department of Mathematics,
University of Rome Tor Vergata, CUP E83C18000100006.
R.L. was partially supported by NSF grant DMS-1800241.
Part of this material is based upon work supported by the National Science Foundation
under Grant No. DMS-1440140 while the authors were in residence at the Mathematical
Sciences Research Institute in Berkeley, California, during the Fall 2018 semester.
}

\thanks{$^*$ Corresponding author: Renato C. Calleja}

\begin{document}
\maketitle

\centerline{\scshape Renato C. Calleja$^*$}
\medskip
{\footnotesize
% please put the address of the first author
 \centerline{Department of Mathematics and Mechanics, IIMAS}
   \centerline{National Autonomous University of Mexico (UNAM), Apdo. Postal 20-72}
   \centerline{C.P. 01000, Mexico D.F. (Mexico)}
} % Do not forget to end the {\footnotesize by the sign }

\medskip

\centerline{\scshape Alessandra Celletti}
\medskip
{\footnotesize
 % please put the address of the second  and third author
 \centerline{Department of Mathematics, University of Rome Tor Vergata}
   \centerline{Via della Ricerca Scientifica 1}
   \centerline{00133 Rome (Italy)}
}

\centerline{\scshape Rafael de la Llave}
\medskip
{\footnotesize
 % please put the address of the second  and third author
 \centerline{School of Mathematics, Georgia Institute of Technology}
   \centerline{686 Cherry St.}
   \centerline{Atlanta GA. 30332-1160 (USA)}
}

\bigskip

% The name of the associate editor will be entered by an editorial staff
% "Communicated by the associate editor name" is not needed for special issue.
 \centerline{(Communicated by the associate editor name)}

%The abstract of your paper

\begin{abstract}
We investigate the existence of whiskered tori in some dissipative systems,
called \sl conformally symplectic \rm systems, having the property that
they transform the symplectic form into a multiple of itself.
We consider a family $f_\mu$ of conformally symplectic maps which
depend on a drift parameter $\mu$.

We fix a Diophantine frequency of the torus and we assume to have a drift
$\mu_0$ and an embedding of the torus $K_0$, which satisfy approximately
the invariance equation $f_{\mu_0} \circ K_0 - K_0 \circ T_\omega$ (where $T_\omega$
denotes the shift by $\omega$). We also assume to have a splitting
of the tangent space at the range of $K_0$ into three bundles.
We assume that the bundles are approximately invariant under $D f_{\mu_0}$ and that
the derivative satisfies some \emph{``rate conditions''}.

Under suitable non-degeneracy conditions, we prove that there exists
$\mu_\infty$, $K_\infty$ and splittings, close to the original ones, invariant
under $f_{\mu_\infty}$. The proof provides an efficient algorithm to construct
whiskered tori. Full details of the statements and proofs are given in \cite{CCL-w}.
\end{abstract}

%\tableofcontents

\section{Introduction}\label{sec:intro}
Whiskered tori for a dynamical system are invariant tori such that the motion
on the torus is conjugated to a rotation and have hyperbolic directions,
exponentially contracting in the future or in the past under the linearized
evolution (\cite{Arnold64,Arnold63a}). Whiskered tori and their invariant manifolds are the key ingredients proposed in
\cite{Arnold64} of the so-called \sl Arnold's diffusion \rm in which solutions of a nearly integrable
system may drift far from their initial values.

Whiskered tori have been widely studied mainly for symplectic systems (see, e.g., \cite{LlaveS}, \cite{FontichLS},
\cite{FontichLSII}); in this paper we go over the results of \cite{CCL-w} and we
consider their existence for \sl conformally symplectic \rm systems (\cite{Banyaga02,CallejaCL13a,DettmannM96,WojtkowskiL98}),
which are characterized by the fact that the symplectic structure is transformed into a
multiple of itself.
Conformally symplectic systems are a very special case of dissipative systems and occur in several
physical examples, e.g. the spin-orbit problem in Celestial Mechanics, Gaussian thermostats,
Euler-Lagrange equations of exponentially discounted systems (\cite{Celletti2010}, \cite{WojtkowskiL98},
\cite{DFIZ2014,DaviniFIZ16}).

The existence of invariant tori in conformally symplectic systems needs an adjustment
of parameters. This leads to consider a family $f_\mu$ of conformally symplectic maps
depending parametrically on $\mu$. Our main result (Theorem~\ref{whiskered}) establishes the
existence of whiskered tori with frequency $\omega$ for $f_\mu$ for some $\mu$;
the Theorem is based on the formulation of an invariance equation
for the parameterization of the torus, say $K=K(\theta)$, for the parameter $\mu$ and
for the splittings of the space. The invariance equation expresses that the parameterization
and the splittings are invariant for the map $f_\mu$. The
main assumption of
Theorem~\ref{whiskered}
is that we are given a sufficiently approximate solution of
\eqref{inv} with an approximately invariant splitting.
We also need to assume that the frequency $\omega$ is Diophantine
and that some non-degeneracy conditions are met.  We note
that the non-degeneracy conditions we need to assume are algebraic
expressions depending only
on the approximate solution and its derivatives. We do not need to
assume any global properties (such as twist) for the whole system.
We also note that the theorem does not make any assumption that the
system considered is close to integrable. Theorems where the main
hypothesis is that there is an approximate solution that have some
condition numbers are called ``a-posteriori'' theorems in the numerical
analysis literature.

The proof of Theorem~\ref{whiskered} is based in showing that a Newton-like
method started on the approximate solution converges.
At each step of the Newton's method, the linearized equation is projected on the hyperbolic and center
subspaces. The equations on the hyperbolic subspaces are solved using a
contraction method (see, e.g., \cite{CCCLwave}).
The invariance equation projected on the center subspace is solved using the so-called
\sl automatic reducibility: \rm taking advantage from the geometry of a conformally symplectic
system, one can introduce a change of coordinates in which the linearized equation
along the center directions can be solved by Fourier methods.

A remarkable result is that
we show that the center bundles of whiskered tori are trivial in the sense of
bundle theory, i.e.
that they are homeomorphic to product bundles. On the other hand,
we allow that the stable and unstable bundles are trivial and there are
examples of this situation.  Note that non-trivial bundles do not seem
to be incorporated in some of the proofs based in normal form theory.

We remark that we do not use transformation
theory as in the pioneering works \cite{Moser67}, \cite{BroerHTB90}, \cite{BroerHS96},
that is we do not perform subsequent changes of variables that transform the system into
a form which admits an invariant torus.

Whiskered tori were studied with a similar approach in \cite{FontichLS}, \cite{FontichLSII};
the results in an
a-posteriori format were proved in \cite{FontichLS} for the case of finite-dimensional
Hamiltonian systems, while generalizations to Hamiltonian lattice systems are
presented in \cite{FontichLSII} and to PDEs in \cite{LlaveS}.

The method introduced in \cite{LlaveGJV05}
(see also \cite{Llave01}, \cite{CallejaCL13a}, and \cite{CanadellH2}
for an application to quasi-periodic normally hyperbolic invariant tori)
has several advantages: it leads to efficient algorithms,
it does not need to work in action-angle variables and it does not assume that the system is
close to integrable. Hence, the approach is suitable to study systems close to breakdown and
in the limit of small dissipation. This allows us to study the analyticity domain of $K$
and $\mu$ as a function of a parameter $\eps$, such that the limit of $\eps$ tending to
zero represents the symplectic case. Note that the limit of dissipation going to zero
is a singular limit. Full dimensional KAM tori in conformally symplectic
systems have also been considered in \cite{Locatelli, LocatelliS15}. The first paper is based on transformation theory and the second includes also
numerical implementations comparing the methods based on transformation
theory and those based on studying \eqref{inv}.

Our second main result, Theorem~\ref{thm:domain}, shows that,
if we introduce an extra perturbative parameter
$\eps$ so that $f_\eps$ is a symplectic map with
a solution $K_0, \mu_0$ of \eqref{inv}, it can be be continued to
$K_\eps$, $\mu_\eps$ which are analytic in
a domain obtained by removing from a ball centered at the origin a sequence of smaller balls,
whose centers lie on a curve and whose radii decrease very fast with their distance from
the origin (see also \cite{CCLdomain}, \cite{Bustamante}). The proof is based on the
construction of Lindstedt series, whose finite order truncation
provides an approximate solution which is used as the approximate solution of the a-posteriori theorem.
We conjecture that such domain is essentially optimal.

The rest of this paper is organized as follows. In Section~\ref{sec:preliminary} we provide
some preliminary notions; Section~\ref{sec:cocycle} presents some properties of cocycles
and invariant bundles; the main result, Theorem~\ref{whiskered}, is stated in Section~\ref{sec:existence};
a sketch of the proof of Theorem~\ref{whiskered} is given in Section~\ref{sec:whiskered};
an algorithm allowing to construct the new approximation
is given in Section~\ref{sec:algorithm};
the analyticity domains
of whiskered tori are presented in Section~\ref{sec:domain}.

\section{Preliminary notions}\label{sec:preliminary}
This section is devoted to introducing the notion of conformally symplectic systems,
the definition of Diophantine vectors, that of invariant rotational tori, and the introduction
of function spaces.

We denote by $\M = \torus^n \times B$ a symplectic manifold of dimension
$2n$ with $B \subseteq \real^n$ an open, simply connected domain with
smooth boundary. We endow $\M$ with the standard
scalar product and a symplectic form $\Omega$, which does not have necessarily the standard form.
In the small dissipation limit (see Section~\ref{sec:domain}), we assume that $\Omega$ is exact.

\begin{definition}\label{defCS}
A diffeomorphism $f:\M\rightarrow \M$ is conformally symplectic, if there exists a function $\lambda$ such that
\begin{equation} \label{conformallysympmap}
f^* \Omega = \lambda\, \Omega\ .
\end{equation}
\end{definition}

We will consider $\lambda$ constant, which is always the case for $n\geq 2$ (\cite{Banyaga02}),
since whiskered tori exist only for $n\geq 2$.

Denoting by $\langle \cdot, \cdot\rangle$ the inner product on $\real^{2n}$, let $J_x$ be the matrix representing $\Omega$ at $x$:
$$
\Omega_x(u, v) = \langle u, J_x v\rangle
$$
with $J^T_x = - J_x$.

Frequency vectors of whiskered tori are assumed to be Diophantine.

\begin{definition}\label{def:DC}
For $\lambda\in\complex$, let $\nu(\lambda;\omega,\tau)$ be defined as
$$
\nu(\lambda;\omega,\tau) \equiv \sup_{k\in\integer^d\setminus\{0\}}
\Big(|e^{2\pi ik\cdot\omega}-\lambda|^{-1}\ |k|^{-\tau}\Big)\ .
$$
We say that $\lambda$ is $\omega$-Diophantine of
class $\tau$ and constant $\nu(\lambda;\omega,\tau)$, if
$$
\nu(\lambda;\omega,\tau)<\infty\ .
$$
\end{definition}
A particular case of the above is when $\lambda = 1$, which corresponds
to the classical definition of $\omega$. In our theorems, we will assume that
$\omega$ is Diophantine and we will consider $\lambda$'s which are Diophantine
with respect to it.

We remark that in Theorem~\ref{whiskered} we will take only $\lambda \in \real$, while in
Theorem~\ref{thm:domain} we will take $\lambda\in\complex$.

To find an invariant torus in a conformally symplectic system, we need to adjust some
parameters (\cite{CallejaCL13a}); hence, we consider a family $f_\mu$ of conformally
symplectic maps depending on a \sl drift \rm parameter $\mu$.

\begin{definition}
Let $f_\mu:\M\rightarrow\M$ be a family of differentiable diffeomorphisms and let
$K: \torus^d \rightarrow \M$ be a differentiable embedding.
Denoting by $T_\omega$ the shift by $\omega\in\real^d$, we say that
$K$ parameterizes an invariant torus for the parameter $\mu$, if the following invariance equation is satisfied:
\begin{equation}\label{inv}
f_\mu \circ K = K \circ T_\omega\ .
\end{equation}
\end{definition}

\medskip

Equation \eqref{inv}, which will be the centerpiece of our study, contains
$K$ and $\mu$ as unknowns; its linearization will be analyzed using a
quasi-Newton method that takes advantage of the geometric properties of conformally symplectic
systems. We remark that if $(K,\mu)$ is a solution, then also $(K\circ T_\alpha,\mu)$ is a solution.
We also show that local uniqueness is obtained by choosing a suitable normalization that fixes $\alpha$.

The analytic function space and a norm is introduced as follows to make estimates
on the quantities involved in the proof.

\begin{definition}\label{def:spaces}
Let $\rho >0$ and let $\torus_\rho^d$ be the set
\[
\torus^d_\rho = \{ z \in \complex^d/\integer^d \, :\ {\rm Re}(z_j)\in\torus\ ,\quad |\Im(z_j)| \leq \rho\ ,\quad j=1,...,d\}\ .
\]
Given a Banach space $X$, let $\A_{\rho}(X)$ be the set of functions
from $\torus^d_\rho$ to $X$, analytic in ${\rm Int}(\torus_\rho^d)$ and
extending continuously to the boundary of $\torus^d_\rho$. We endow $\A_{\rho}$ with the following norm,
which makes $\A_\rho$ a Banach space:
$$
\|f\|_{\A_{\rho}} = \sup_{z\in\torus_\rho^d}\ |f(z)|\ .
$$
\end{definition}

The norm of a vector valued function $g=(g_1,\ldots,g_n)$ is defined as
$\|g\|_{\A_{\rho}}=\sqrt{\|g_1\|^2_{\A_{\rho}}+\ldots+\|g_n\|^2_{\A_{\rho}}}$,
while the norm of an $n_1\times n_2$ matrix valued function $G$ is defined as
$\|G\|_{\A_{\rho}}=\sup_{\chi\in\real^{n_2}_+,|\chi|=1}
\sqrt{\sum_{i=1}^{n_1}(\sum_{j=1}^{n_2}\|G_{ij}\|_{\A_{\rho}}\, \chi_j)^2}$.

\section{Cocycles and invariant bundles}
\label{sec:cocycle}

Given an approximate solution of \equ{inv}, we will be led
to reduce the error and hence to study products of the form
\begin{equation} \label{product}
\Gamma^j \equiv f_\mu \circ K \circ  T_{(j-1)\omega} \times  \cdots \times  Df_\mu\circ K\ ,
\end{equation}
which are quasi-periodic cocycles of the form
\begin{equation}\label{cocyclerotation}
\Gamma^j = \gamma_\theta\circ T_{(j-1) \omega} \times \cdots \times \gamma_\theta\
\end{equation}
with $\gamma_\th = Df_\mu\circ K(\th)$.
The cocycle \eqref{cocyclerotation} satisfies the property:
$\Gamma^{j+m} = \Gamma^j\circ T_{m\omega}\ \Gamma^m$.
The study of the invariance equation strongly depends on the asymptotic growth
of the cocycle \equ{product}, which leads to the following definition (\cite{SackerS74,Coppel78}).

\begin{definition} \label{def:trichotomy}
The cocycle \eqref{product} admits an exponential trichotomy if there exists a decomposition
\begin{equation}\label{splitting}
\real^n =  E_{\th}^s \oplus E_{\th}^c \oplus E_{\th}^u\ ,
\qquad \theta\in\torus^d\ ,
\end{equation}
rates of decay
$\lambda_- < \lambda_c^- \le \lambda_c^+ <  \lambda_+$,
$\lambda_- < 1 < \lambda_+$
and a constant $C_0 > 0$, such that
\begin{equation}\label{growthrates}
\begin{split}
& v\in  E_{\th}^s\ \iff   |\Gamma^j(\th) v | \le C_0 \lambda_-^j |v|,\quad j \ge 0 \\
& v\in  E_{\th}^u\ \iff   |\Gamma^j(\th) v | \le C_0 \lambda_+^j |v|,\quad j \le 0 \\
&v\in  E_{\th}^c\ \iff
\begin{matrix}
|\Gamma^j(\th) v | \le C_0 (\lambda_c^-)^j  |v|, \quad j \ge 0 \\
\ \ |\Gamma^j(\th) v | \le  C_0 (\lambda_c^+)^j |v|, \quad  j \le  0\ . \\
\end{matrix}
\end{split}
\end{equation}
\end{definition}

Given a splitting as in \eqref{splitting}, we denote by
$\Pi^s(\th), \Pi^c(\th), \Pi^u(\th)$ the projections,
depending on the whole splitting, on
$E_\th^s$, $E_\th^c$, $E_\th^u$. Let us now consider two
nearby splittings, $E$, $\tE$; then, for each space
in $\tE$, we can find a linear function
$A^\sigma_\th: E^\sigma_\th \rightarrow E^{\hat\sigma}_\th$
(where $E^{\hat\sigma}_\th$ is the sum of the spaces
in the splitting not indexed by $\sigma$), such that
\begin{equation} \label{graph}
\tE^\sigma_\th = \{ v \in \real^n, v = x + A^\sigma_\th x \, | \,
x \in E^\sigma_\th \}\ .
\end{equation}
% We will take as reference splitting the initial approximate splitting
% satisfying the hypothesis of Theorem~\ref{whiskered}.

Denoting by $P^\perp_{E_\th}$ the orthogonal projections, the distance between $E$ and $\tE$ is defined as
$$
{\rm dist}_\rho(E, \tE)  = \| P^\perp_{E_\th} - P^\perp_{\tE_\th} \|_{\A_\rho}\ .
$$

From \eqref{growthrates} it is possible to show (\cite{HirschPS77}) that the splittings depend
continuously (H\"older) on $\th$; bootstrapping the regularity,
the splittings are analytic. Therefore, the projections $\Pi^\sigma$,
$\sigma =  s,u,c$, are  uniformly bounded (\cite{SackerS74}).
Also, we remark that the bundles characterized by
\eqref{growthrates} are invariant:
$\gamma_\th E^\sigma_\th = E^{\sigma}_{\th + \omega}$ (\cite{CCL-w}).

\subsection{Approximately invariant splittings}
For a splitting $E^s_\th  \oplus E^u_\th \oplus E^c_\th$ and
a cocycle $\gamma_\th$, let $\gamma^{\sigma, \sigma'}_\th$ be
\beq{sigmasigma}
\gamma^{\sigma, \sigma'}_\th =
 \Pi^\sigma_{\th +\omega}  \gamma_\th \Pi^{\sigma'}_\th\ ;
\eeq
hence, the splitting is invariant under the cocycle if and only if
\[
\gamma_\th^{\sigma, \sigma'} \equiv 0\ , \qquad \sigma \ne \sigma'\ .
\]
The lack of invariance of the splitting under the cocyle $\gamma$ is measured by the quantity
$$
\I_\rho(\gamma, E)  \equiv  \max_{\sigma, \sigma' \in \{s,c,u\} \atop \sigma \ne \sigma'}
\sup_{\torus^d_\rho} \| \gamma_\th^{\sigma, \sigma'} \|_{\rho}\ .
$$
Now, we introduce a notion of hyperbolicity for approximately invariant splittings.

\begin{definition}\label{approximatehyperbolic}
Let $\gamma$ be a cocycle and $E$ an approximately invariant splitting.
Then, $\gamma$ is approximately  hyperbolic w.r.t. $E$, if the cocycle
\[
\tilde \gamma_\theta =
\begin{pmatrix}
& \gamma^{s,s}_\theta & 0 & 0 \\
&0 &  \gamma^{c,c}_\theta & 0  \\
&0 & 0 & \gamma^{u,u}_\theta
\end{pmatrix}
\]
satisfies \equ{growthrates} with $\gamma^{\sigma, \sigma}$ as in \equ{sigmasigma}.
\end{definition}

The following Lemma~\ref{lem:closing} shows that if we have an approximately invariant
splitting for an approximately hyperbolic cocycle, then there exists a true invariant
splitting.

\begin{lemma}\label{lem:closing}
Fix an analytic reference splitting on $\torus^d_\rho$ and let $\U$ be a
sufficiently small neighborhood of this splitting, so
that all the splittings can be parameterized as graphs of
linear
maps $A^\sigma_\th$ as in \eqref{graph}
with $\|A^\sigma_\th\|_\rho  < M_1$.

Let $E$ be an analytic splitting  in the neighborhood $\U$.

Let $\gamma$ be an analytic cocycle over a rotation defined on
$\torus^d_\rho$ with $\| \gamma \|_\rho  < M_2$ for $M_2\in\real_+$.

Assume that $E$ is approximately invariant under $\gamma$:
$$
\I_\rho(\gamma, E) \le \eta
$$
and that $\gamma$ is approximately hyperbolic for
the reference splitting as in Definition~\ref{approximatehyperbolic}.

Then, there is a locally unique splitting $\tilde E$ close to $E$, invariant under
$\gamma$, which satisfies the trichotomy of Definition~\ref{def:trichotomy}, and
such that
$$
\dist_\rho( E, \tilde E) \le C \eta\ ,
$$
where $C$, $\eta$ can be chosen uniformly and depending only on $M_1$, $M_2$.
\end{lemma}

We refer to \cite{CCL-w} for the proof of the closing Lemma~\ref{lem:closing}, which is based on the standard method of
writing the new spaces as the graphs of linear  maps
$A_x^\sigma: E^\sigma  \rightarrow E_x^{\hat \sigma}$
(were $E^{\hat \sigma}_x$ denotes the sum of the spaces in the splitting that are
different from $E^\sigma_x$).
The fact that these spaces are invariant can be transformed into fixed point
equations that can be solved by the contraction mapping principle.
We refer to \cite{CCL-w} for details.

\section{Existence of whiskered tori}\label{sec:existence}

Whiskered tori are defined as follows.

\begin{definition}
\label{def:whiskered}
Let $f_\mu:\M\rightarrow\M$ be a family of conformally symplectic maps with
conformal factor $\lambda$.
We say that $K: \torus^d \rightarrow \M$ represents a whiskered torus when
for some $\omega\in\real^d$:
\begin{enumerate}
\item
$K$ is the embedding of a rotational torus: $f_\mu \circ K = K \circ T_\omega$.
\item
The cocycle $Df_\mu\circ K$ over the rotation $T_\omega$
admits a trichotomy as in \equ{growthrates} with rates
$\lambda_-, \lambda_c^-, \lambda_c^+, \lambda_+$.
\item
The rates satisfy
$\lambda_c^- \le \lambda \le \lambda_c^+$.
\item
The spaces $E^c_\th$ in \equ{splitting} have dimension $2d$.
\end{enumerate}
\end{definition}

Theorem~\ref{whiskered} below states
the existence of whiskered tori by solving the invariance equation \equ{inv}.

Let $K$, $\mu$ be an approximate solution of \eqref{inv} with a small error term $e$:
$f_{\mu}\circ K - K\circ T_\omega = e$. Let $\Delta$, $\beta$ be some corrections, so
that $K'=K+\Delta$, $\mu'=\mu+\beta$ satisfy the invariance equation with an error
quadratically smaller. This is obtained, provided $\Delta$, $\beta$ satisfy
$$
(Df_{\mu} \circ K)\ \Delta - \Delta \circ T_\omega
+ (D_\mu f_{\mu}) \circ K \beta  = -e\ .
$$

\begin{theorem}\label{whiskered}
Let $\omega \in \D_d(\nu, \tau)$, $d\leq n$, be as in \eqref{def:DC},
let $f_\mu:\M\rightarrow\M$, $\mu \in \real^d$, be a family of real analytic, conformally symplectic mappings
as in \equ{conformallysympmap} with $0<\lambda<1$. We make the following assumptions.

$(H1)$ Appproximate solution:

Let $(K_0,\mu_0)$ with $K_0 :\torus^d \to \M$, $K_0 \in\A_\rho$,
and $\mu_0 \in \real^d$ represent an approximate whiskered torus for $f_{\mu_0}$ with frequency $\omega$:
$$
\|f_{\mu_0}\circ K_0  - K_0 \circ T_\omega\|_{\A_\rho}\leq \E\ ,\qquad \E>0\ .
$$

To ensure that the composition $f_\mu\circ K$ can be defined, we assume that there exists a domain
$\U \subset \complex^n/\integer^n \times \complex^n $ such that
for all $\mu$ with $| \mu -\mu_0| \le \eta$, $f_\mu$ has domain $\U$ and
\begin{equation}
\label{compositions}
\dist( K_0(\torus^d_\rho), \complex^n/\integer^n \times \complex^n
\setminus \U) \ge \eta.
\end{equation}

$(H2)$
Approximate splitting:

For all $\th\in\torus^d_\rho$, there exists a splitting of the tangent
space of the phase space, depending analytically on the angle
$\th\in\torus^d_\rho$; the bundles are approximately invariant under the
cocycle $\gamma_\th = Df_{\mu_0} \circ K_0(\th)$, i.e.
$\I_\rho(\gamma,E)\leq\E_h$, $\E_h>0$.

$(H3)$ Spectral condition for the bundles (exponential trichotomy):

For all  $\th\in\torus_\rho^d$ the spaces in $(H2)$
are approximately hyperbolic for the cocycle
$\gamma_\th$.

$(H3')$ Since we are dealing with conformally symplectic
systems, we assume:
$$
\lambda_- < \lambda \lambda_+< \lambda_c^-\ ,\qquad \lambda_c^-\leq \lambda\leq\lambda_c^+\ .
$$

$(H4)$ The dimension of the center subspace is $2d$.

$(H5)$ Non--degeneracy:

Let $N(\th)=(DK(\th)^T DK(\th))^{-1}$, $P(\th)=DK(\th) N(\th)$, $\chi(\th)=DK(\th)^T(J_c)^{-1}\circ K(\th) DK(\th)$,
and let
\beq{esse1}
S(\th)=P(\th+\omega)^T Df_\mu \circ K(\th) (J^c)^{-1}\circ K(\th)P(\th)
- N(\th+\omega)^T \chi(\th + \omega) N(\th+\omega)\ \lambda\,\Id_d\ .
\eeq
Let $M$ be defined as
\beq{M}
M(\th) = [ DK(\th)\ |\  (J^c)^{-1}\circ K(\th)\ DK(\th) N(\th)]\ .
\eeq

We assume that an explicit $d \times d$ matrix $\S$, formed
by algebraic operations (and solving cohomology equations)
from the derivatives of
the approximate solution, is invertible.

Let $\alpha=\alpha(\tau)$ be an explicit number and assume that for some $0< \delta < \rho$, we have
$\E \le \delta^{2 \alpha} \E^*$, $\E_h  \le \E^*_h$, where $\E$, $\E^*$ depend on
$\nu$, $\tau$, $C_0$, $\lambda_+$, $\lambda_-$, $\lambda_c^+$, $\lambda_c^-$,
$\| \Pi^{s/u/c}_{\th}\|_{\A_\rho}$, $\| DK_0\|_{\A_\rho}$, $\| (DK_0^{T} DK_0)^{-1}\|_{\A_\rho}$

Then, there exists an exact solution $(K_e,\mu_e)$, which satisfies
$$
f_{\mu_e}\circ K_e-K_e\circ T_\omega=0
$$
with
\beqano
\|K_e-K_0\|_{\A_{\rho-2\delta}}\leq C\E\delta^{-\tau} ,\qquad |\mu_e-\mu_0|\leq
C \E\ .
\eeqano
Moreover, the invariant torus $K_e$ is hyperbolic, namely there exists an invariant
splitting
\[
\T_{K_e(\th)}\M = E_{\th}^s \oplus
 E_{\th}^c \oplus E_{\th}^u\ ,
\]
satisfying Definition~\ref{def:trichotomy} and which is close to
the original one:
\[
\| \Pi^{s/u/c}_{K_e(\th) } -
\Pi^{s/u/c}_{K_0(\th)}   \|_{\A_{\rho - 2 \delta}} \le C (\E\delta^{-\tau} + \E_h)\ .
\]
Finally, the hyperbolicity constants associated to the invariant splitting
of the invariant torus (denoted by a tilde) are close to those of the approximately
invariant splitting of the approximately invariant torus
(see (H1), (H2)):
$$
| \lambda_\pm - \widetilde\lambda_\pm| \le C ( \E\delta^{-\tau}  + \E_h)\ , \quad
| \lambda^\pm_c  - \widetilde\lambda^{\pm}_c | \le C ( \E \delta^{-\tau}  + \E_h)\ .
$$
\end{theorem}

\begin{remark}
Some consequences of the geometry are the following (see \cite{CCL-w}).

$\bullet$ The stable/unstable exponential rates given by
the set of Lyapunov multipliers $\{\lambda_i\}_{i=1}^{2d}$
satisfy the pairing rule
$$
\lambda_i\ \lambda_{i+d}=\lambda\ .
$$

$\bullet$ Invariant tori satisfy the
isotropic property: the symplectic form restricted to the invariant torus is zero.

$\bullet$ Because of the conformally symplectic structure,
the symplectic form is non-degenerate
when restricted to the center bundle $E^c_{K(\th)}$.
\end{remark}

An important result is that the bundle $E_{K(\th)}^c$ near a rotational invariant torus satisfying
our hypotheses (notably that the dimension of the fibers of the bundle
is $2d$) is trivial, that is, the bundle is isomorphic to
a product bundle. Precisely, we can show that if $K$ is an approximate solution
of \equ{inv}, we can find a linear operator $\B_\theta:Range(DK(\theta))\rightarrow E_{K(\theta)}^c$,
such that $E_{K(\theta)}^c$ is the range under $\Id+\B_\theta$ of the tangent bundle of the torus.

\section{A sketch of the proof of Theorem~\ref{whiskered}}\label{sec:whiskered}

We now proceed to sketch the proof of Theorem~\ref{whiskered} (see Section~\ref{sec:sketch}), which uses
the so-called automatic reducibility presented in Section~\ref{sec:automatic}.
The proof leads to the algorithm described in Section~\ref{sec:algorithm}.
We refer to \cite{CCL-w} for full details.

\subsection{Automatic reducibility}\label{sec:automatic}

We assume that there exists an invariant splitting of
the tangent space of $\M$ at $K(\th)$, $\T_{K(\th)}\M$ with $\th\in\torus^d$:
$$
\T_{K(\th)} \M=E_{K(\th)}^s \oplus E_{K(\th)}^c \oplus E_{K(\th)}^u\ .
$$
Taking the derivative of \equ{inv} we get
\beq{derinv}
Df_\mu\circ K(\th)\, DK(\th) - DK\circ T_\omega(\th) =0\ ,
\eeq
which shows that $Range(DK(\th))\subset E_{K(\th)}^c$ and hence:
\begin{equation}
\label{isotropiccenter}
DK^T(\th) J^c\circ K(\th)  DK(\th) = 0\ ,
\end{equation}
where $J^c$ is the $2n\times 2n$ matrix of the embeddings of the center space into the ambient space.
Due to \eqref{isotropiccenter}, the dimension of the range of $M$ in \equ{M} is $2d$ and,
from $(H4)$, we have:
\begin{equation}\label{rangeM}
\Range(M(\th)) = E^c_{K(\th)}\ .
\end{equation}
Hence, there exists a matrix $\B(\th)$ such that
\beq{red}
Df_\mu \circ K(\th) M(\th) = M(\th +\omega)\ \B(\th)\ ,
\eeq
where $\B(\th)$ is upper triangular with constant matrices on the diagonal.
From \equ{derinv}, the first columm of
$\B$ is $\begin{pmatrix} \Id_d  \\  0 \end{pmatrix}$.
From \eqref{rangeM}, setting $v(\th) = (J^c)^{-1}\circ K(\th)\ DK(\th) N(\th)$,
we have
\begin{equation}\label{secondcolumn}
Df_\mu\circ K(\th) \, v(\th) = DK(\th+\omega) S(\th) +  v(\th + \omega) U(\th)\ ,
\end{equation}
where $U=U(\th)$ is obtained multiplying \eqref{secondcolumn} on the right by
$DK^T(\th + \omega) J^c\circ K(\th)$ and using \eqref{isotropiccenter}:
\beq{UU}
U(\th)=DK^T(\th + \omega ) J^c\circ K(\th +\omega)Df_\mu\circ K(\th) \, v(\th)\ .
\eeq
From the conformally symplectic and invariance properties of
the center foliation, we obtain:
\[
Df^T_\mu(x) J_{f(x)}^c Df_\mu(x)  = \lambda J_{f(x)}^c\ ,
\]
from which $J_{f(x)}^c  Df_\mu(x) (J_{x}^c)^{-1}  = \lambda Df_\mu^{-T}(x)$.

Hence, we see that the the left hand side of \eqref{UU} is equal to $\lambda$,
thus showing that
$$
U(\th) = \lambda\ .
$$
Defining $S$ as in \equ{esse1}, we can write \equ{red} as
\beq{red2}
Df_\mu \circ K(\th) M(\th) = M(\th +\omega)
\begin{pmatrix} \Id_d & S(\th)\\ 0&\lambda\Id_d \end{pmatrix}\ .
\eeq

\subsection{Sketch of the proof}\label{sec:sketch}

Once the automatic reducibility leading to \equ{red2} is established, we can proceed
to sketch the proof of Theorem~\ref{whiskered}.

We start with an approximate solution of the invariance
equation which is approximately
hyperbolic  and look for a correction to $K$ and $\mu$, such that the error
of the invariance of the new embedding and the new parameter
is, roughly, the square of the original error in a smaller domain; this is the content of
the following Proposition.

\begin{proposition}\label{pro:iterative}
Let $f_\mu:\M\rightarrow\M$, $\mu\in\real^d$, $d\leq n$,
be a family of real-analytic, conformally symplectic maps as in Theorem~\ref{whiskered} with
$0<\lambda<1$. Let $\omega\in D_d(\nu,\tau)$.

Let $(K,\mu)$, $K :\torus^d \to \M$, $K\in\A_\rho$, be an approximate solution, such that
\beq{invw}
f_\mu\circ K(\th)-K\circ T_\omega(\th)=e(\th)
\eeq
and let $\E=\|e\|_{\A_\rho}$.

Let $E^{s/c/u}_{K(\th)}$ be an approximately invariant hyperbolic splitting based on $K$,
such that $\I_\rho(\gamma,E^{s/c/u}_{K(\th)})<\E_h$.
Assume that $(K,\mu)$ satisfy assumptions $(H2)$-$(H3)$-$(H3')$-$(H4)$-$(H5)$ of Theorem~\ref{whiskered}
and that $\E$, $\E_h$ are sufficiently small.

Then, there exists an exact invariant splitting $\widetilde E^{s/c/u}_{K(\th)}$ with associated
cocycle $\widetilde \gamma_\th^{\sigma,\sigma}$, such that
$$
dist_\rho(E^{s/c/u}_{K(\th)},\widetilde E^{s/c/u}_{K(\th)})\leq C\E_h\ ,\qquad
\|\gamma_\th^{\sigma,\sigma}-\widetilde \gamma_\th^{\sigma,\sigma}\|_{\A_\rho}\leq C\E_h\ .
$$
Furthermore, we can find corrections $\Delta$, $\beta$, such that $K'=K+\Delta$, $\mu'=\mu+\beta$ satisfy
$$
f_{\mu'}\circ K'(\th)-K'\circ T_\omega(\th)=e'(\th)
$$
with
$$
\|e'\|_{\A_{\rho-\delta}}\leq C\ \delta^{-2\tau}\ \E^2\ ,\qquad
\|\Delta\|_{\A_{\rho-\delta}}\leq C\ \delta^{-\tau}\ \E\ ,\qquad
|\beta|\leq C\E\ .
$$
Moreover, the splitting $\widetilde E^{s/c/u}_{K(\th)}$ is approximately invariant for $Df_{\mu'}\circ K'$.
\end{proposition}

The proof of Proposition~\ref{pro:iterative} is based on the following ideas.
Expanding in Taylor series the invariance equation for $K'$, $\mu'$, we have:
\beqano
f_{\mu'}\circ K'(\th)-K'(\th+\omega)&=&f_\mu\circ K(\th)+Df_\mu\circ K(\th)\,\Delta(\th)+D_\mu f_\mu \circ K(\th)\,\beta\nonumber\\
&-&K(\th+\omega)-\Delta(\th+\omega)+O(\|\Delta\|^2)+O(|\beta|^2)\ .
\eeqano
Using \equ{invw}, the new error is quadratically smaller if the corrections $\Delta$, $\beta$ satisfy
\beq{EW}
Df_\mu\circ K(\th)\ \Delta(\th)+D_\mu f_\mu\circ K(\th)\ \beta-\Delta(\th+\omega)=-e(\th)\ .
\eeq
The solution of \equ{EW} is obtained by projecting it on the hyperbolic
and center spaces, and
using the invariant splitting \equ{splitting}. Let $K_e$ be the exact
solution of \equ{invw};
we assume that the cocycle $Df_\mu\circ K_e$ admits an invariant splitting as in
\equ{splitting}. For the initial step, this follows from $(H2)$ and the closing Lemma~\ref{lem:closing}, while in
subsequent steps, we observe that the exactly invariant splitting
for one step will be approximately invariant for the corrected one,
so that we can apply again Lemma~\ref{lem:closing} to
restore the invariance.

Denoting by $\Delta^{\xi}(\th)\equiv \Pi_{K(\th+\omega)}^{\xi}\Delta(\th)$,
$e^{\xi}(\th)\equiv \Pi_{K(\th+\omega)}^{\xi}e(\th)$ with $\xi=s,c,u$, we have
\beq{WB}
Df_\mu\circ K(\th)\ \Delta^{\xi}(\th)+
\Pi^{\xi}_{K(\th+\omega)} D_\mu f_\mu\circ K(\th)\beta-\Delta^{\xi}(\th+\omega)
=-e^{\xi}(\th)\ ,
\eeq
which contains $\Delta^s$, $\Delta^c$, $\Delta^u$, $\beta$ as unknowns.
The equation for $\Delta^c$ allows to determine $\Delta^c$ and $\beta$. In fact, from
$$
\Delta^c=M\ W^c\ ,
$$
recalling \equ{red}, the approximate solution satisfies \equ{red2} up to an error term, say $R=R(\th)$:
\beq{redW}
Df_\mu\circ K(\th)\ M(\th)=M(\th+\omega)\B(\th)+R(\th)
\eeq
with
\[
\|R\|_{\A_{\rho-\delta}}\leq C \delta^{-1}\ \|e\|_{\A_\rho}\ .
\]
Using \equ{WB} and \equ{redW}, one obtains:
\beq{CW}
\left(%
\begin{array}{cc}
  \Id_d & S(\th) \\
  0 & \lambda\Id_d \\
 \end{array}%
\right)W^c(\th)-W^c\circ T_\omega(\th)=-\widetilde e^c(\th)-\widetilde A^c(\th)\beta\ ,
\eeq
where $\widetilde e^c(\th)\equiv M^{-1}\circ T_\omega(\th) e^c(\th)$,
$\widetilde A^c(\th)\equiv M^{-1}\circ T_\omega(\th)\
\Pi_{K(\th+\omega)}^c D_\mu f_\mu\circ K(\th)$.
Next, we define $\widetilde A^c\equiv [\widetilde A^c_1|\widetilde A^c_2]$,
$\overline{W^c}$ is the average of $W^c$, $(W^c)^0\equiv W^c-\overline{W^c}$ and,
being $\noaverage{W_2^c}$ an affine function of $\beta$, we let
$\noaverage{W_2^c}=\noaverage{W_a^c}+\beta\noaverage{W_b^c}$ for some
functions $W_a^c$, $W_b^c$. With this setting, \equ{CW} becomes
\beqa{eqW}
\noaverage{W_1^c}(\th) - \noaverage{W_1^c}\circ T_\omega(\th)
 &=& - \noaverage{S W_2^c}(\th)
 -\noaverage { \widetilde e_1^c}(\th) - \noaverage{\widetilde A_1^c}(\th) \beta\nonumber\\
\lambda \noaverage{W_a^c}(\th)  -  \noaverage{W_a^c} \circ T_\omega(\th)
&=& - \noaverage{\widetilde e_2^c}(\th)\nonumber\\
\lambda \noaverage{W_b^c}(\th)  -  \noaverage{W_b^c} \circ T_\omega(\th)
&=& - \noaverage{\widetilde A_2^c}(\th)\ ,
\eeqa
whose solution for $\noaverage{W_1^c}$, $\noaverage{W_a^c}$, $\noaverage{W_b^c}$ is found
by using standard results (see, e.g., \cite{CallejaCL13a}).

Taking the average of \equ{eqW}, recollecting the last two equations in a single equation
for $\noaverage{W_2^c}$, leads to solve the following system
\beq{sys}
\left(\begin{array}{cc}
  \overline{ S} & \overline{S (W_b^c)^0}+\overline{\widetilde A_1^c} \\
  (\lambda -1)\Id_d & \overline{\widetilde A_2^c} \\
\end{array}\right)
\left(\begin{array}{c}
  \overline{W_2^c} \\
  \beta \\
\end{array}\right)=
\left(\begin{array}{c}
  -\overline{ S (W_a^c)^0   }- \overline{\widetilde e_1^c} \\
  -\overline{\widetilde {e_2^c}} \\
\end{array}\right)\ .
\eeq
Using the non-degeneracy condition (H5),%\eqref{non-degeneracyW},
allows to find a solution of \equ{sys} and, hence,
to determine $\overline{W_2^c}$, $\beta$.

Next, we solve \equ{WB} for the stable subspace. Denoting by $\th'=T_\omega(\th)$,
$\widetilde e^{s}(\th')\equiv \Pi_{K(\th')}^{s} e\circ T_{-\omega}(\th')$, equation \equ{WB} becomes
$$
Df_\mu(K\circ T_{-\omega}(\th'))\ \Delta^{s}(T_{-\omega}(\th'))+
\Pi_{K(\th+\omega)}^s D_\mu f_\mu(K\circ T_{-\omega}(\th'))\beta
-\Delta^{s}(\th')=-\widetilde e^{s}(\th')\ ,
$$
which can be solved for $\Delta^s$ in the form
\beqano
&&\Delta^s(\th')=\widetilde e^s(\th')+\sum_{k=1}^\infty \Big(Df_\mu(K\circ T_{-\omega}(\th'))\times \cdots \times
Df_\mu(K\circ T_{-k\omega}(\th'))\Big)\ \widetilde e^s(T_{-k\omega}(\th'))\nonumber\\
&+&\Pi_{K(\th+\omega)}^s D_\mu f_\mu(K\circ T_{-\omega}(\th'))\beta\nonumber\\
&+&\sum_{k=1}^\infty \Big(D f_\mu(K\circ T_{-\omega}(\th'))\times...\times
D f_\mu(K\circ T_{-k\omega}(\th'))\
\Pi_{K(\th+\omega)}^s D_\mu f_\mu(K\circ T_{-(k+1)\omega}(\th')\Big)\ \beta\ ,
\eeqano
where the series in the last term converges in $\A_\rho$, due to the growth rates
\equ{growthrates}.

In a similar way, one can solve the equation for the unstable subspace, thus obtaining
\beqano
&&\Delta^u(\th)=-\sum_{k=0}^\infty \Big((Df_\mu)^{-1}(K(\th))\times...\times
(Df_\mu)^{-1}(K\circ T_{k\omega}(\th))\Big)\  e^u(T_{k \omega}(\th))\nonumber\\
&-&\sum_{k=0}^\infty \Big((Df_\mu)^{-1}(K(\th))\times...\times
(Df_\mu)^{-1}(K\circ T_{k\omega}(\th))\
\Pi_{K(\th+\omega)}^u D_\mu f_\mu(K\circ T_{k\omega}(\th)\Big)\beta\ .
\eeqano

Simple estimates lead to state that $\S$ is (H5), the norm of the projections,
the change in the rates and the constant $C_0$ in \equ{growthrates} slightly change after one iterative step;
denoting by $\widehat \gamma^{\sigma, \sigma}_\th$ the cocycle associated to $K'$, $\mu'$, one has:
\beqano
\|\S'\|_{\A_{\rho-\delta}}&\leq&  \|\S\|_{\A_{\rho}}+C\delta^{-\tau} \|e\|_{\A_{\rho}}\nonumber\\
\|\Pi^{s/c/u}_{K'(\th)}-\Pi^{s/c/u}_{K(\th)}\|_{\A_{\rho}}&\leq& C\|K'-K\|_{\A_{\rho}}
\leq C\delta^{-\tau}\|e\|_{\A_{\rho}}\nonumber\\
\|\widehat \gamma^{\sigma, \sigma}_\th-\widetilde \gamma^{\sigma, \sigma}_\th\|_{\A_{\rho-\delta}}&\leq&
C(\delta^{-\tau}\E+\E_h)\ .
\eeqano
The last issue to prove Theorem~\ref{whiskered} is to show that the inductive step can be iterated
infinitely many times and that it converges to
the true solution, provided the initial error is sufficiently small. This is a standard KAM
argument, which is proved by introducing a sequence $\{K_h,\mu_h\}$ of approximate solutions
on shrinking domains and imposing a smallness condition on the size of the
initial error $\|e\|_{\A_{\rho}}$.

\section{The algorithm for the new approximation}\label{sec:algorithm}
The proof of Theorem~\ref{whiskered} leads to the following algorithm, which allows to
construct the improved approximation, given $f_\mu$, $\omega$, $K_0$, $\mu_0$.
We fix an integer $L_0$, which denotes the maximum number of terms which are computed in
the infinite series defining $\Delta^s$ and $\Delta^u$.
Each step is denoted as $a \leftarrow b$, meaning that the quantity $a$ is determined from $b$.  Note that the number of steps is less than 40 and that all
the steps involve just calling a standard function, so that the coding is
sort of straightforward.

\begin{algorithm}
Let $f_\mu$, $\omega$, $K_0$, $\mu_0$ be as in the previous sections and let $L_0\in\integer$:

$\bullet$ $\chi_1$ $\leftarrow$  $f_{\mu_0}\circ K_0$

$\bullet$ $\chi_2$ $\leftarrow$ $K_0\circ T_\omega$

$\bullet$ $e$ $\leftarrow$ $\chi_1-\chi_2$

$\bullet$ $e^{s/c/u}$ $\leftarrow$ $\Pi_{\th+\omega}^{s/c/u}e$

$\bullet$ $\gamma$ $\leftarrow$ $Df_{\mu_0}\circ K_0$

$\bullet$ $\widetilde\gamma$ $\leftarrow$ $D_\mu f_{\mu_0}\circ K_0$

$\bullet$ $\alpha$ $\leftarrow$ $DK_0$

$\bullet$ $N$ $\leftarrow$ $[\alpha^T\alpha]^{-1}$

$\bullet$ $\widetilde J$ $\leftarrow$ $(J^c)^{-1}\circ K$

$\bullet$ $M$ $\leftarrow$ $[\alpha|\widetilde J\ \alpha\ N]$

$\bullet$ $\widetilde M$ $\leftarrow$ $M^{-1}\circ T_\omega$

$\bullet$ $\widetilde e^c$ $\leftarrow$ $\widetilde M\ e^c$

$\bullet$ $P$ $\leftarrow$ $\alpha N$

$\bullet$ $\chi$ $\leftarrow$ $\alpha^T\ \widetilde J\ \alpha$

$\bullet$ $\Lambda$ $\leftarrow$ $\lambda\ \Id_d$

$\bullet$ $S$ $\leftarrow$ $(P\circ T_\omega)^T\gamma\widetilde J\ P-(N\circ T_\omega)^T(\chi\circ T_\omega)^T\ N\circ T_\omega\ \Lambda$

$\bullet$ $\widetilde A^c$ $\leftarrow$ $\widetilde M\ \Pi_{\th+\omega}^c\ \widetilde\gamma$

$\bullet$ $(W_a^c)^o$ solves $\lambda(W_a^c)^o-(W_a^c)^o\circ T_\omega=-(\widetilde e_2^c)^o$

$\bullet$ $(W_b^c)^o$ solves $\lambda(W_b^c)^o-(W_b^c)^o\circ T_\omega=-(\widetilde A_2^c)^o$

$\bullet$ Find $\overline{W_2^c}$, $\beta$ solving
\beqano
\overline{S}\ \overline{W_2^c}+(\overline{S(W_b^c)^o}+\overline{\widetilde A_1^c})\beta&=&-\overline{S(W_a^c)^o}
-\overline{\widetilde e_1^c}\nonumber\\
(\lambda-1)\overline{W_2^c}+\overline{\widetilde A_2^c}\beta&=&-\overline{\widetilde e_2^c}
\eeqano

$\bullet$ $(W_2^c)^o$ $\leftarrow$ $(W_a^c)^o+\beta(W_b^c)^o$

$\bullet$ $W_2^c$ $\leftarrow$ $(W_2^c)^o+\overline{W_2^c}$

$\bullet$ $(W_1^c)^o$ solves $(W_1^c)^o-(W_1^c)^o\circ T_\omega=-(S\ W_2^c)^o-(\widetilde e_1^c)^o-(\widetilde A_1^c)^o\beta$

$\bullet$ $\Delta^c$ $\leftarrow$ $M^c\ W^c$

$\bullet$ $\mu_1$ $\leftarrow$ $\mu_0+\beta$

$\bullet$ Compute $\widetilde\Gamma_k=\gamma^{-1}\times ... \times\gamma^{-1}\circ T_{k\omega}$ for $k=0,...,L_0$

$\bullet$ Compute $e_k^u=e^u\circ T_{k\omega}$ for $k=0,...,L_0$

$\bullet$ Compute $\widetilde \gamma_k=\Pi_{\th+\omega}^u\widetilde\gamma\circ T_{k\omega}$ for $k=0,...,L_0$

$\bullet$ $\Delta^u$ $\leftarrow$ $-\sum_{k=0}^{L_0}(\widetilde \Gamma_k e_k^u+\widetilde \Gamma_k\widetilde \gamma_k\beta)$

$\bullet$ $\theta'$ $\leftarrow$ $T_\omega(\theta)$

$\bullet$ Compute $\Gamma_k=\gamma\circ T_{-\omega}(\th')\times ...\times \gamma\circ T_{-k\omega}(\th')$ for $k=1,...,L_0$

$\bullet$ Compute $\widetilde e_k^s=\widetilde e^s\circ T_{-k\omega}(\th')$ for $k=1,...,L_0$

$\bullet$ Compute $\widetilde\gamma'=\Pi_{\th'}^s\widetilde\gamma\circ T_{-\omega}(\th')$ for $k=1,...,L_0$

$\bullet$ Compute $\widehat\gamma_k'=\Pi_{\th'}^s\widetilde\gamma\circ T_{-(k+1)\omega}(\th')$ for $k=1,...,L_0$

$\bullet$ $\Delta^{s '}$ $\leftarrow$ $\widetilde e^s(\th')+\sum_{k=1}^{L_0} \Gamma_k\widetilde e_k^s+\widetilde \gamma'\beta
+\sum_{k=1}^{L_0}\Gamma_k\widehat\gamma_k'\beta$

$\bullet$ $\Delta^s$ $\leftarrow$ $\Delta^{s '}\circ T_{-\omega}$

$\bullet$ $K_1'$ $\leftarrow$ $K_0+\Delta^c+\Delta^u+\Delta^s$.

\end{algorithm}

\section{Domains of analyticity and Lindstedt expansions of
whiskered tori}\label{sec:domain}

The study of domains of analyticity of
whiskered tori of conformally symplectic systems in the limit of
small dissipation is similar to  that developed in
\cite{CCLdomain}, but adding the hyperbolicity.
The main idea is to compute an asymptotic expansion (Lindstedt series), which can be used
as starting point for the application of
Theorem~\ref{whiskered}.

The Lindstedt series expansions to order $N$ of $K$, $\mu$, $A^\sigma$
satisfy the invariance equation up to
an error bounded by $C_N |\eps|^{N+1}$. Then, we apply
Theorem~\ref{whiskered} for $\eps$ belonging to a domain with good
Diophantine properties of $\lambda$. Hence, we are able to prove
that there exists a true solution $K$, $\mu$ and that
$$
\|K^{[\le N]}- K\|, |\mu^{[ \le N]}-\mu | \le \tilde C_N |\eps|^{N+1}
$$
in the domain, thus showing that the Lindstedt
series are asymptotic expansions of the true solution.
The quantities $K^{[\le N]}$, $\mu^{[ \le N]}$ denote the truncations to
order $N$ in $\eps$ (see \equ{Hlam} below) of the Lindstedt series expansions.

Let $f_{\mu_\varepsilon,\varepsilon}:\M\rightarrow\M$ be a family of maps, such that
$$
f_{\mu_\varepsilon,\varepsilon}^*\,\Omega=\lambda(\varepsilon)\Omega\ ,
$$
where the conformal factor $\lambda$ is taken as
\beq{Hlam}
\lambda(\eps) = 1 + \alpha \eps^a + O(|\eps|^{a+1})
\eeq
for some $a>0$ integer and $\alpha\in\complex\setminus\{0\}$.

Recalling Definition~\ref{def:DC}, we introduce the following sets, where the Diophantine constants
allow to start an iterative convergent procedure (see \cite{CCLdomain}).

\begin{definition}\label{def:sets}
For $A>0$, $N\in\integer_+$, $\omega\in\real^d$, let the set $\G=\G(A;\omega,\tau,N)$ be defined as
$$
\G(A;\omega,\tau,N)\equiv\{\eps\in\complex:\ \nu(\lambda(\eps);\omega,\tau)\ |\lambda(\eps)-1|^{N+1}\leq A\}\ .
$$
For $r_0\in\real$, let
\beq{Gr0}
\G_{r_0}(A;\omega,\tau,N)=\G\cap \{\eps\in\complex:\  |\eps| \le r_0\}\ .
\eeq
\end{definition}

\begin{figure}[h]
\centering
\includegraphics[width=6truecm,height=6truecm]{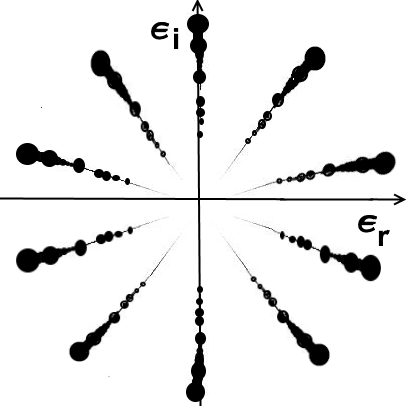}
\caption{A representation of the set $\G$ given by the complement of the black circles,
whose radii have been rescaled for graphical reasons.
We took $d=1$, $\tau=1$, $a=5$.}
\label{figure2}
\end{figure}

We prove that $K$ and $\mu$
are analytic in a domain $\G_{r_0}$ as in \equ{Gr0} for a
sufficiently small $r_0$. This domain is obtained by removing from a
ball centered at zero a sequence of smaller balls whose centers lie
along smooth lines going through the origin (see Figure~\ref{figure2}).
The removed balls have radii decreasing faster than any power of the distance
of their center from the origin.
Like in \cite{CCLdomain}, we conjecture that this domain is essentially optimal.

\begin{theorem}\label{thm:domain}
Let $f_{\mu,\eps}:\M\rightarrow\M$ with
$\mu\in\Gamma$ with $\Gamma\subseteq \complex^d$ open, $d\leq n$, $\eps\in\complex$, be a family of conformally
symplectic maps with conformal factor satisfying \equ{Hlam} with
$\alpha \in \real$, $\alpha \ne 0$, $a \in \nat$. Let $\omega \in \D_d(\nu, \tau)$.

(A1) Assume that for $\mu=\mu_0,\eps= 0$ the map $f_{\mu_0, 0}$ admits a whiskered invariant torus, namely

(A1.1) there exists an embedding $K_0:\torus^d \rightarrow \M$, $K_0  \in \A_\rho$ for some
$\rho>0$, such that
$$
f_{\mu_0,0} \circ K_0 = K_0 \circ T_\omega\ ;
$$

(A1.2) there exists a splitting $T_K(\th)\M =
E^s_{K(\th)} \oplus E^c_{K (\th)}  \oplus E^u_{K(\th)}$, which is invariant under the cocycle
$\gamma^0_\th = D f_{\mu_0,0} \circ K_0(\th) $
and satisfies Definition~\ref{def:trichotomy}.
The ratings of the  splitting  satisfy the assumptions (H3), (H3')
and (H4) of Theorem~\ref{whiskered}.

(A.2) The function $f_{\mu, \eps}(x)$ is
analytic in all of its arguments and that the analyticity domains
are large enough, namely:

(A2.1) both $K_0(\th)$ and the splittings $E^{s,c,u}_{K(\th)}$
considered as a function of $\th$ are  in $\A_{\rho_0}$ for some $\rho_0 > 0$;

(A2.2) there is a domain $\U\subset\complex^n/\integer^n \times \complex^n $ such
that for $|\eps| \le \eps^*$ and all $\mu$ with
$| \mu - \mu_0| \le \mu^*$ the function $f_{\mu, \eps}$ is
defined in $\U$ and we also have \eqref{compositions}.

(A3) The non-degeneracy condition (H5) of Theorem~\ref{whiskered} is satisfied by
the invariant torus.

\noindent
Then:

\medskip

B.1) We can compute formal power series expansions
\[
K_\eps^{[\infty]}= \sum_{j=0}^\infty\eps^j  K_j
\qquad
\mu_\eps^\infty =  \sum_{j=0}^\infty\eps^j  \mu_j
\]
satisfying \eqref{inv} and such that for any $0 < \rho' \le \rho$ and $N\in\nat$, we have
$$
|| f_{\mu_{\eps}^{[\le N]}, \eps} \circ K^{[\le N]}_\eps -  K^{[\le N]}_\eps \circ T_\omega ||_{\rho'} \le C_N |\eps|^{N+1}\ .
$$

B.2) We can compute the formal power series expansions
\[
A^{\sigma,\infty}_\eps= \sum_{j=0}^\infty\eps^j  A_j^\sigma\ ,\qquad
A^{\sigma}_j(\th) : E^\sigma_0(\th) \rightarrow \E^{\hat \sigma}_0(\th)\ ,\qquad
\sigma = \st, \hat \st, \un, \hat \un
\]
with $A^\sigma_j \in \A_\rho$ satisfying the equations
%\eqref{invariance2-s}, \eqref{invariance2-u}
for invariant dichotomies  in the sense of power series.

\medskip

B.3) For the set $\G_{r_0}$ as in
\equ{Gr0} with $r_0$ sufficiently small and for $0 < \rho' < \rho$,
there exists $K_\varepsilon:\G_{r_0}\rightarrow\A_{\rho'}$, $\mu_\varepsilon:\G_{r_0}\rightarrow\complex^d$,
analytic in the interior of
$\G_{r_0}$ taking values in $\A_{\rho'}$,
extending continuously to the boundary of $\G_{r_0}$ and
such that for $\eps\in\G_{r_0}$ the invariance equation is satisfied exactly:
$$
f_{\mu_\eps, \eps} \circ K_\eps - K_\eps \circ T_\omega=0\ .
$$
Furthermore, the exact solution admits the formal series in A) as an asymptotic expansion, namely for
$0<\rho'<\rho$, $N \in \nat$, one has:
$$
||K^{[\le N]}_\eps -  K_\eps||_{\rho'}  \le C_N |\eps|^{N+1} \ ,\qquad
|\mu^{[\le N]}_\eps -  \mu_\eps|  \le C_N |\eps|^{N+1}\ .
$$
\end{theorem}

We refer to \cite{CCL-w} for the proof of Theorem~\ref{thm:domain}.
Here we limit to give a graphical description as in Figure~\ref{figure2}
of the set $\G$, which is the complement of the black circles
with centers along smooth lines going through the origin and with radii
decreasing very fast as the centers go to zero.

\begin{comment}
\appendix

\section{The solution of the cohomological equation}\label{sec:app}
A cohomological equation of the form
\begin{equation}
w(\varphi + \omega) -\lambda w(\varphi) = \eta (\varphi)
\label{difference1}
\end{equation}
for some functions $w$ and $\eta$ with $\eta$ having zero average, can be solved through
the following result, which is standard in KAM theory (see, e.g., \cite{CallejaCL13a} for the proof).

\begin{lemma}\label{neutral}
Consider \eqref{difference1} for $\lambda \in [A_0, A_0^{-1}]$ and for some
$0 < A_0< 1$; let $\omega\in \D_d(\nu,\tau)$. Assume that $\eta \in
\A_{\rho}$, $\rho>0$.
Then, there is one and  only one solution of \eqref{difference1}
with zero average and if $\varphi \in \A_{\rho-\delta}$ for some
$0<\delta<\rho$, then
$$
\|\varphi\|_{\A_{\rho-\delta}} \le C\, \nu\, \delta^{-\tau} \|\eta\|_{\A_{\rho}}\ ,
$$
for some constant $C$ depending on $A_0$, $d$ and which is uniform in $\lambda$ and independent on the
Diophantine constant $\nu$.
\end{lemma}

\end{comment}

\def\cprime{$'$} \def\cprime{$'$} \def\cprime{$'$} \def\cprime{$'$}

%\bibliographystyle{alpha}
%\bibliography{cclwhiskers}

\begin{thebibliography}{CCCdlL17}

\bibitem[{A}rn63]{Arnold63a}
V.~I. {A}rnol'd.
\newblock Proof of a theorem of {A}. {N}. {K}olmogorov on the invariance of
  quasi-periodic motions under small perturbations.
\newblock {\em Russian Math. Surveys}, 18(5):9--36, 1963.

\bibitem[Arn64]{Arnold64}
V.I. Arnold.
\newblock Instability of dynamical systems with several degrees of freedom.
\newblock {\em Sov. Math. Doklady}, 5:581--585, 1964.

\bibitem[Ban02]{Banyaga02}
A.~Banyaga.
\newblock Some properties of locally conformal symplectic structures.
\newblock {\em Comment. Math. Helv.}, 77(2):383--398, 2002.

\bibitem[BC18]{Bustamante}
A.~P. Bustamante and R.~Calleja.
\newblock Computation of domains of analyticity for the dissipative standard
  map in the limit of small dissipation.
\newblock {\em Preprint}, 2018.

\bibitem[BHS96]{BroerHS96}
H.~W. Broer, G.~B. Huitema, and M.~B. Sevryuk.
\newblock {\em Quasi-Periodic Motions in Families of Dynamical Systems. {\rm
  Order Amidst Chaos}}.
\newblock Springer-Verlag, Berlin, 1996.

\bibitem[BHTB90]{BroerHTB90}
H.~W. Broer, G.~B. Huitema, F.~Takens, and B.~L.~J. Braaksma.
\newblock Unfoldings and bifurcations of quasi-periodic tori.
\newblock {\em Mem. Amer. Math. Soc.}, 83(421):viii+175, 1990.

\bibitem[CCCdlL17]{CCCLwave}
Renato~C. Calleja, Alessandra Celletti, Livia Corsi, and Rafael de~la Llave.
\newblock Response solutions for quasi-periodically forced, dissipative wave
  equations.
\newblock {\em SIAM J. Math. Anal.}, 49(4):3161--3207, 2017.

\bibitem[CCdlL13]{CallejaCL13a}
Renato~C. Calleja, Alessandra Celletti, and Rafael de~la Llave.
\newblock A {KAM} theory for conformally symplectic systems: efficient
  algorithms and their validation.
\newblock {\em J. Differential Equations}, 255(5):978--1049, 2013.

\bibitem[CCdlL17]{CCLdomain}
Renato~C Calleja, Alessandra Celletti, and Rafael de~la Llave.
\newblock Domains of analyticity and {L}indstedt expansions of {K}{A}{M} tori
  in some dissipative perturbations of {H}amiltonian systems.
\newblock {\em Nonlinearity}, 30(8):3151, 2017.

\bibitem[CCdlL18]{CCL-w}
Renato~C. Calleja, Alessandra Celletti, and Rafael de~la Llave.
\newblock Existence of whiskered {K}{A}{M} tori of conformally symplectic
  systems.
\newblock {\em Preprint,  https://arxiv.org/abs/1901.07483}, 2019.

\bibitem[Cel10]{Celletti2010}
Alessandra Celletti.
\newblock {\em Stability and Chaos in Celestial Mechanics}.
\newblock Springer-Verlag, Berlin; published in association with Praxis
  Publishing Ltd., Chichester, 2010.

\bibitem[CH17]{CanadellH2}
Marta Canadell and \`Alex Haro.
\newblock Computation of quasiperiodic normally hyperbolic invariant tori:
  rigorous results.
\newblock {\em J. Nonlinear Sci.}, 27(6):1869--1904, 2017.

\bibitem[Cop78]{Coppel78}
W.~A. Coppel.
\newblock {\em Dichotomies in stability theory}.
\newblock Lecture Notes in Mathematics, Vol. 629. Springer-Verlag, Berlin-New
  York, 1978.

\bibitem[DFIZ16a]{DFIZ2014}
Andrea Davini, Albert Fathi, Renato Iturriaga, and Maxime Zavidovique.
\newblock Convergence of the solutions of the discounted equation: the discrete
  case.
\newblock {\em Math. Z.}, 284(3-4):1021--1034, 2016.

\bibitem[DFIZ16b]{DaviniFIZ16}
Andrea Davini, Albert Fathi, Renato Iturriaga, and Maxime Zavidovique.
\newblock Convergence of the solutions of the discounted {H}amilton-{J}acobi
  equation: convergence of the discounted solutions.
\newblock {\em Invent. Math.}, 206(1):29--55, 2016.

\bibitem[dlL01]{Llave01}
R.~de~la Llave.
\newblock A tutorial on {KAM} theory.
\newblock In {\em Smooth ergodic theory and its applications ({S}eattle, {WA},
  1999)}, volume~69 of {\em Proc. Sympos. Pure Math.}, pages 175--292. Amer.
  Math. Soc., Providence, RI, 2001.

\bibitem[dlLGJV05]{LlaveGJV05}
R.~de~la Llave, A.~Gonz{\'a}lez, {\`A}.~Jorba, and J.~Villanueva.
\newblock K{AM} theory without action-angle variables.
\newblock {\em Nonlinearity}, 18(2):855--895, 2005.

\bibitem[dlLS18]{LlaveS}
R.~de~la Llave and Y.~Sire.
\newblock An a posteriori {K}{A}{M} theorem for whiskered tori in {H}amiltonian
  partial differential equations with applications to some ill-posed equations.
\newblock {\em Arch. Rational. Mech. Anal.,
  https://doi.org/10.1007/s00205-018-1293-6}, 2018.

\bibitem[DM96]{DettmannM96}
C.~P. Dettmann and G.~P. Morriss.
\newblock Proof of {L}yapunov exponent pairing for systems at constant kinetic
  energy.
\newblock {\em Phys. Rev. E}, 53(6):R5545--R5548, 1996.

\bibitem[FdlLS09]{FontichLS}
Ernest Fontich, Rafael de~la Llave, and Yannick Sire.
\newblock Construction of invariant whiskered tori by a parameterization
  method. {P}art {I}: maps and flows in finite dimensions.
\newblock {\em J. Differential Equations,}, 246:3136--3213, 2009.

\bibitem[FdlLS15]{FontichLSII}
Ernest Fontich, Rafael de~la Llave, and Yannick Sire.
\newblock Construction of invariant whiskered tori by a parameterization
  method. {P}art {II}: {Q}uasi-periodic and almost periodic breathers in
  coupled map lattices.
\newblock {\em J. Differential Equations}, 259(6):2180--2279, 2015.

\bibitem[HPS77]{HirschPS77}
M.W. Hirsch, C.C. Pugh, and M.~Shub.
\newblock {\em Invariant manifolds}.
\newblock Springer-Verlag, Berlin, 1977.
\newblock Lecture Notes in Mathematics, Vol. 583.

\bibitem[LS15]{LocatelliS15}
Ugo Locatelli and Letizia Stefanelli.
\newblock Quasi-periodic motions in a special class of dynamical equations with
  dissipative effects: a pair of detection methods.
\newblock {\em Discrete Contin. Dyn. Syst. Ser. B}, 20(4):1155--1187, 2015.

\bibitem[Mos67]{Moser67}
J.~Moser.
\newblock Convergent series expansions for quasi-periodic motions.
\newblock {\em Math. Ann.}, 169:136--176, 1967.

\bibitem[SL12]{Locatelli}
Letizia Stefanelli and Ugo Locatelli.
\newblock {K}olmogorov's normal form for equations of motion with dissipative
  effects.
\newblock {\em Discrete Contin. Dynam. Systems}, 17(7):2561--2593, 2012.

\bibitem[SS74]{SackerS74}
Robert~J. Sacker and George~R. Sell.
\newblock Existence of dichotomies and invariant splittings for linear
  differential systems. {I}.
\newblock {\em J. Differential Equations}, 15:429--458, 1974.

\bibitem[WL98]{WojtkowskiL98}
M.~P. Wojtkowski and C.~Liverani.
\newblock Conformally symplectic dynamics and symmetry of the {L}yapunov
  spectrum.
\newblock {\em Comm. Math. Phys.}, 194(1):47--60, 1998.

\end{thebibliography}

\end{document}